\documentclass[11pt]{amsart}
\usepackage{geometry}                
\usepackage[parfill]{parskip}    
\usepackage{graphicx}
\usepackage{amssymb}
\usepackage{epstopdf}
\newtheorem{theorem}{Theorem}[section]

\theoremstyle{definition}
\newtheorem{definition}[theorem]{Definition}

\newcommand{\R}{\mathbb{R}}
\newcommand{\N}{\mathbb{N}}

\title[Heteroclinic limit cycle]
      {Global attraction and repulsion of a heteroclinic limit cycle in three dimensional Kolmogorov maps}

\author[Zhanyuan Hou]{}

\subjclass{Primary: 37C29; Secondary: 37C65, 37C70, 37E30}
 \keywords{competitive maps, carrying simplex, heteroclinic limit cycle, global attraction and repulsion}

 \email{z.hou@londonmet.ac.uk}

\begin{document}

\maketitle

\centerline{\scshape Zhanyuan Hou }
\medskip
{\footnotesize
 \centerline{School of Computing and Digital Media, London Metropolitan University,}
   \centerline{166-220 Holloway Road, London N7~8DB, UK}
} 

\medskip


\begin{abstract}
There is a recent development in the carrying simplex theory for competitive maps: under some weaker conditions a map has a modified carrying simplex (one of the author's latest publications). In this paper, as one of the applications of the modified carrying simplex theory, a criterion is established for a three dimensional Kolmogorov map to have a globally repelling (attracting) heteroclinic limit cycle. As a concrete example, a discrete competitive model is investigated to illustrate the above criteria for global repulsion (attraction) of a hetericlinic limit cycle.
\end{abstract}

\textbf{Note.} This paper has been accepted for publication in Proceedings of ICDEA 2021.


\section{Introduction}\label{Sec1}
We are concerned with the global asymptotic behaviour of the discrete dynamical system
\begin{equation}\label{e1}
	x(n) = T^n(x), \quad x\in C,\quad n\in \N,
\end{equation}
where $C = \R^N_+ = [0, +\infty)^N$, $\N = \{0, 1, 2, \ldots\}$ and $T: C \to C$ is the Kolmogorov map having the form
\begin{equation}\label{e2}
	T_i(x) = x_if_i(x), \quad i\in I_N = \{1, 2, \ldots, N\}
\end{equation}
and $f\in C^1(C, C)$ with $f_i(x)>0$ for all $x\in C$ and $i\in I_N$. It is known that system (\ref{e1}) with (\ref{e2}) is a typical mathematical model for the population dynamics of a community of $N$ species, where each $x_i(n)$ represents the population size or density at time $n$ (at the end of $n$th time period), and the function $f_i(x)$ denotes the per capita growth rate, of the $i$th species. If $\frac{\partial f_i}{\partial x_j} \leq 0$ for all $i, j\in I_N$ with $i\not= j$, then  (\ref{e1}) models the population dynamics of a community of competitive species. 

In this paper, we focus on a recent development on the carrying simplex theory for competitive maps. The main aim is to provide a criterion for a three dimensional system to have a globally attracting (repelling) heteroclinic limit cycle based on the existence of a modified carrying simplex. This demonstrates that the new development in the carrying simplex theory has the potential to have various applications for a broader class of systems. 

The paper is organised as follows: section 2 is for the recent development on carrying simplex theory, section 3 presents a criterion for global attraction (repulsion) of a  heteroclinic limit cycle, section 4 deals with a Ricker model as a concrete example illustrating the criterion given in section 3, section 5 is for conclusion.

\section{Carrying simplex of competitive Kolmogorov maps}\label{Sec2}
Research on system (\ref{e1}) with (\ref{e2}) and its various particular instances as models has been flourishing  in the last few decades. The carrying simplex theory and its various applications is one of the important and influential developments. This theory was originally established by Hirsch \cite{Hir1} (see \cite{Hou1} and \cite{Hou2} for latest update) for competitive Kolmogorov systems of differential equations. Since then the idea of a carrying simplex for discrete systems gradually appeared in literature (see \cite{WaJi1}, \cite{WaJi2}, \cite{JiMiWa}, \cite{DiWaYa} for example). But a more accepted theorem for existence and uniqueness of a carrying simplex was given by Hirsch \cite{Hir2} without proof. Then Ruiz-Herrera \cite{Her} presented a more general theorem covering Hirsch's result with a complete proof.

\begin{definition}
	A nonempty set $\Sigma \subset C$ is called a carrying simplex of (\ref{e1}) if it is a compact invariant hypersurface homeomorphic to $\Delta^{N-1} =\{x\in C: x_1 + \cdots +x_N =1\}$ such that every trajectory except the origin is asymptotic with a trajectory in $\Sigma$. 
\end{definition}

If (\ref{e1}) has a carrying simplex $\Sigma$, then, since it attracts all the points $x\in C\setminus\{0\}$ so that the limit set $\omega(x) \subset \Sigma$, the asymptotic dynamics of (\ref{e1}) on $C$ is essentially described by the dynamics on $\Sigma$. Many applications of the carrying simplex theory actually utilise this attracting feature of $\Sigma$. 

For any $x, y\in C$, we write $x\leq y$ or $y \geq x$ if $x_i\leq y_i$ for all $i\in I_N$; $x<y$ or $y>x$ if $x\leq y$ but $x\not= y$; $x\ll y$ or $y\gg x$ if $x_i < y_i$ for all $i\in I_N$. 

\begin{definition}
	The map $T$ given by (\ref{e2}) is said to be retrotone (or competitive) in a subset $X\subset C$ if for any $x, y\in X$, $T(x) < T(y)$ implies $x_i < y_i$ for all $i\in I(y) = \{j\in I_N: y_j\not= 0\}$.	
\end{definition}

Let $[0, r] = \{x\in C: 0\leq x\leq r\}$. The theorem below is Theorem 6.1 in \cite{Her}.

\begin{theorem}\label{The1.1}
	Assume that $T$ with $T([0, r]) \subset [0, r]$ for some $r\gg 0$ satisfies the following conditions:
	\begin{itemize}
		\item[(i)] For each $i\in I_N$, the map $T$ restricted to the positive half $x_i$-axis has a fixed point $q_ie_i$ with $q_i>0$, $e_i$ the $i$th standard unit vector and $q\ll r$.
		\item[(ii)] $T$ is retrotone and locally one to one in $[0, r]$.
		\item[(iii)] For any $x, y\in [0, r]$, if $T(x) < T(y)$ then, for each $j\in I_N$, either $x_j=0$ or $f_j(x) >f_j(y)$.
	\end{itemize}
	Then the map admits a carrying simplex $\Sigma$.
\end{theorem}

Theorem \ref{The1.1} can be only applied to (\ref{e1}) restricted to the space $[0, r]\subset C$ if no condition for $T$ on $C\setminus[0, r]$ is provided. However, if for any compact set $S\subset C$ there is a $k\in\N$ such that $T^k(S) \subset [0, r]$, then Theorem \ref{The1.1} can be applied directly to the system on $C$.

When $f$ on $C$ is a $C^1$ map, $T$ is also a $C^1$ map with Jacobian matrix
\begin{equation}\label{e3}
	DT(x) = \textup{diag}(f_1(x), \ldots, f_N(x))(I - M(x)),
\end{equation}
where $I$ is the identity matrix and 
\begin{equation}\label{e4}
	M(x)= (M_{ij}(x))= \left(-\frac{x_i}{f_i(x)}\frac{\partial f_i}{\partial x_j}(x)\right)_{N\times N}.
\end{equation}
Then, by Lemma 4.1, Corollary 6.1 and Remark 6.4 in \cite{Her}, Theorem \ref{The1.1} has the following version with easily checkable conditions.

\begin{theorem}\label{The1.2}
	Assume that $T$ satisfies the following conditions:
	\begin{itemize}
		\item[(i)] For each $i\in I_N$, the map $T$ restricted to the positive half $x_i$-axis has a fixed point $q_ie_i$ with $q_i>0$, $e_i$ the $i$th standard unit vector and $q\ll r$ for some $r\in C$.
		\item[(ii)] All entries of the Jacobian $Df$ are negative.
		\item[(iii)] The spectral radius of $M(x)$ satisfies $\rho(M(x)) <1$ for all $x\in [0, q]\setminus\{0\}$.
	\end{itemize}
	Then the map admits a carrying simplex $\Sigma$.
\end{theorem}

A more user-friendly variation of Theorem \ref{The1.2} given by Jiang and Niu \cite[Theorem 3.1]{JiNi2} is the above theorem with $\rho(M(x)) <1$ replaced by an inequality involving $\frac{\partial f_i}{\partial x_j}(x)$. (See (\ref{e7}) or (\ref{e8}) below.) 

We note that condition (ii) in Theorem \ref{The1.2} is very restrictive; it excludes the possibility of applying the theorem to systems with some zero entries of $Df$. But actually, condition (ii) is too strong and unnecessary, a compact invariant set attracting all the points of $C\setminus\{0\}$ with most of the features of a carrying simplex may still exist even if $\frac{\partial f_i}{\partial x_j} =0$ for some distinct $i, j\in I_N$. The author's recent work \cite{Hou3} provides a weaker sufficient condition for (\ref{e1}) to have a modified carrying simplex. Let $\dot{C}$ denote the interior of $C$ and let $B$ be either $C$ or a positively invariant $[0, r]$ for some $r\in\dot{C}$. 

\begin{definition}
	A nonempty set $\Sigma\subset B\setminus\{0\}$ is called a modified carrying simplex of (\ref{e1}) if $\Sigma$ meets the following requirements.
	\begin{itemize}
		\item[(i)] $\Sigma$ is compact, invariant and homeomorphic to $\Delta^{N-1}$ by radial projection.
		\item[(ii)] $\Sigma$ attracts all the points of $B\setminus\{0\}$, i.e. $\omega(x)\subset \Sigma$ for each $x\in B\setminus\{0\}$.  
	\end{itemize}
	Moreover, if $x$ is below $\Sigma$ with a nonempty support $I(x)\subset I_N$, then there is a $y\in\Sigma$ with $I(y)=I(x)$ such that $\lim_{n\to +\infty}(T^n(x)-T^n(y)) =0$.
\end{definition}

Note that the main difference between modified carrying simplex and the carrying simplex in literature is that the latter requires every trajectory in $B\setminus\{0\}$ to be asymptotic to one in $\Sigma$ whereas the former requires every nontrivial trajectory below $\Sigma$ to be asymptotic to one in $\Sigma$ and $\Sigma$ to attract all the points of $B\setminus\{0\}$. Obviously, the concept of a modified carrying simplex is less restrictive and it includes carrying simplex as a particular class.

\begin{definition}
	The map $T: C\to C$ defined by (\ref{e2}) is said to be weakly retrotone (or weakly competitive) in a subset $X\subset C$ if for $x, y\in X$ with $T(x) >T(y)$ and $I = \{i\in I_N: T_i(x)>T_i(y)\}$, then $x>y$ and $x_i>y_i$ for all $i\in I$.
\end{definition}

Comparing this with the definition of retrotone we see that if $T$ is retrotone then it is weakly retrotone, but not vice versa. The theorem below is Theorem 2.3 in \cite{Hou3}.

\begin{theorem}\label{The2.1}
	Assume that $T$ defined by (\ref{e2}) with $T([0, r]) \subset [0, r]$ for some $r\in \dot{C}$ satisfies the following conditions:
	\begin{itemize}
		\item[(i)] For each $i\in I_N$, the map $T$ restricted to the positive half $x_i$-axis has a fixed point $q_ie_i$ with $q_i>0$ and $q\ll r$.
		\item[(ii)] $T$ is weakly retrotone and locally one to one in $[0, r]$.
		\item[(iii)]  For any $x, y\in [0, r]$, if $T(x) < T(y)$ with $I = \{i\in I_N: T_i(x)<T_i(y)\}$ then, for each $j\in I$, either $x_j=0$ or $f_j(x) >f_j(y)$.
	\end{itemize}
	Then 0 is a repellor with the basin of repulsion $\mathcal{B}(0)\subset [0, r]$, (\ref{e1}) has a unique modified carrying simplex $\Sigma$ and $\Sigma = \overline{\mathcal{B}(0)} \setminus (\{0\}\cup\mathcal{B}(0))$. Moreover, for each $p\in\Sigma$ and every $q\in [0, r]\setminus\{0\}$ with $q<p$, we have $\alpha(q)\subset \pi_i$ provided $q_i< p_i$.
\end{theorem}

Now utilising $DT$ and $Df$, we obtain conditions which guarantee conditions (ii) and (iii) and the following version of Theorem \ref{The2.1} with easily checkable conditions. Consider the matrix $M(x)$ given by (\ref{e4}) and
\begin{equation}
	\tilde{M}(x)= (\tilde{M}_{ij}(x))= \left(-\frac{x_j}{f_i(x)}\frac{\partial f_i}{\partial x_j}(x)\right)_{N\times N}. \label{E4}
\end{equation}
The theorem below is Theorem 2.4 in \cite{Hou3}. Note that (\ref{E4}) was used by some other authors in literature (e.g. \cite{JiNi2}).

\begin{theorem}\label{The2.2}
	Assume that $T$ given by (\ref{e2}) satisfies the following conditions:
	\begin{itemize}
		\item[(i)] For each $i\in I_N$, the map $T$ restricted to the positive half $x_i$-axis has a fixed point $q_ie_i$ with $q_i>0$ and $q\ll r$ for some $r\in \dot{C}$.
		\item[(ii)] The entries of the Jacobian $Df$ satisfy
		\begin{equation}\label{e6}
			\forall x\in [0, r], \forall i, j\in I_N,\; \frac{\partial f_i}{\partial x_j}(x) \leq 0,
		\end{equation}
		and $f_i$ is strictly decreasing in $x_i\in [0, r_i]$ for $x\in [0, r]$.
		\item[(iii)] For each $x\in [0, q]\setminus\{0\}$, either $\rho(M(x))<1$ for $M(x)$ given by (\ref{e4}) or $\rho(\tilde{M}(x))<1$ for $\tilde{M}(x)$ given by (\ref{E4}).
	\end{itemize}
	Then 0 is a repellor with the basin of repulsion $\mathcal{B}(0)\subset [0, r]$, (\ref{e1}) has a unique modified carrying simplex $\Sigma$ and $\Sigma = \overline{\mathcal{B}(0)} \setminus (\{0\}\cup\mathcal{B}(0))$. Moreover, for each $p\in\Sigma$ and every $q\in [0, r]\setminus\{0\}$ with $q<p$, we have $\alpha(q)\subset \pi_i$ provided $q_i< p_i$.
\end{theorem}

Note that
\begin{equation}\label{e7}
	f_i(x) + \sum_{j=1}^Nx_i\frac{\partial f_i}{\partial x_j}(x)>0, \quad \forall i\in I_N,
\end{equation}
implies $\rho(M(x))<1$ and
\begin{equation}\label{e8}
	f_i(x) + \sum_{j=1}^Nx_j\frac{\partial f_i}{\partial x_j}(x)>0, \quad \forall i\in I_N,
\end{equation}
implies $\rho(\tilde{M}(x)) <1$ (see \cite{Hou3}). Then condition (iii) in Theorem \ref{The2.2} is satisfied if either (\ref{e7}) or (\ref{e8}) holds.

Theorem \ref{The1.2} is a popular result in discrete competitive systems and has a large number of applications due to the important and interesting features of a carrying simplex: compact, invariant, unordered ($p\leq q$ implies $p=q$ for $p, q\in \Sigma$), homeomorphic to $\Delta^{N-1}$ by radial projection, and attracting all the points of  $C\setminus\{0\}$. The following are just a few examples. Ruiz-Herrera \cite{Her} investigated exclusion and dominance utilizing the existence of a carrying simplex. Jiang and Niu \cite{JiNi1}, \cite{JiNi2} and Gyllenberg et al. \cite{GyJiNiYa1}, \cite{GyJiNiYa2} dealt with some well known three-dimensional competitive models. Based on the existence of a carrying simplex, they classified the systems into 33 topologically equivalent classes and gave a phase portrait on $\Sigma$ for each class. Jiang, Niu and Wang \cite{JiNiWa} studied heteroclinic cycles via carrying simplex. Balreira et al. \cite{BaElLu} and Gyllenberg et al. \cite{GyJiNi3} provided criteria for global stability of an interior fixed point based on the existence of a carrying simplex. Baigent \cite{Bai1}, \cite{Bai2} investigated the geometric feature of a carrying simplex and found conditions for $\Sigma$ to be convex. Baigent and Hou \cite{BaHo} and Hou \cite{Hou} provided split Lyapunov function method and geometric method for global stability. Although these methods were not based on the existence of a carrying simplex, comments and comparisons with those using carrying simplex were made there.

With the introduction of a modified carrying simplex, we expect that Theorem \ref{The2.2} can be applied to a broader class of systems in various applications. In section 3 of \cite{Hou3}, geometric criteria for dominance and vanishing species was provided based on the existence of a modified carrying simple. In the next section, as an application of Theorem \ref{The2.2} we shall prove a theorem for (\ref{e1}) to have a globally attracting (repelling) heteroclinic limit cycle.

\section{Global attraction and repulsion of a heteroclinic limit cycle}
In this section, we are going to find sufficient conditions for (\ref{e1}) to have a globally attracting (repelling) heteroclinic limit cycle.  We assume that the conditions of Theorem \ref{The2.2} are satisfied so that (\ref{e1}) has a unique modified carrying simplex $\Sigma$.

By a {\em heteroclinic cycle} we mean a closed curve that is topologically a circle consisting of fixed points $p_{(i)}$ for $i\in I_N$, together with heteroclinic curves  $\ell_i$ connecting $p_{(i)}$ to $p_{(i+1)}$ (here $p_{(N+1)}=p_{(1)}$). By a {\em heteroclinic limit cycle} $\Gamma$ we mean a heteroclinic cycle $\Gamma$ with an attracting (or repelling) neighbourhood $N(\Gamma)$ (restricted to $\dot{C}$) such that $\omega(x^0)=\Gamma$ (or $\alpha(x^0) = \Gamma$) for all $x^0\in N(\Gamma)$. We first define the concept of {\em globally} attracting or repelling heteroclinic limit cycle.

\begin{definition} 
	We say that a heteroclinic cycle $\Gamma_0$ of (\ref{e1}) is a
	\begin{itemize}
		\item locally attracting (repelling) heteroclinic limit cycle if there is a neighbourhood $V\subset \dot{C}$ ($V\subset\textup{int}\Sigma$) of $\Gamma_0$ such that $\omega(x^0) =\Gamma_0$ $(\alpha(x^0) =\Gamma_0)$ for all $x^0\in V$;
		\item globally attracting (repelling) heteroclinic limit cycle if $\omega(x^0) = \Gamma_0$ $(\alpha(x^0) =\Gamma_0)$
		for all $x^0\in \dot{C} \setminus W^s(p)$ $(x^0\in\textup{int}\Sigma \setminus \{p\})$, where $W^s(p)$ is the one-dimensional stable manifold of a fixed point $p\in \dot{C}$.
	\end{itemize}
\end{definition}

For the case where $N=3$ and system (\ref{e1}) admits a carrying simplex $\Sigma$, the three axial fixed points $Q_1, Q_2, Q_3$ are the only fixed points on $\partial\Sigma$, and $\partial\Sigma$ is a heteroclinic cycle, Jiang, Niu and Wang \cite{JiNiWa} gave a simple sufficient condition $\det(\theta) < 0 (>0)$ for $\partial\Sigma$ to be locally attracting (repelling), where $\theta =(\ln f_i(Q_j))_{3\times 3}$. The main issue we address here is when the heteroclinic cycle is globally attracting (repelling) limit cycle.

For convenience, let
\[
\pi_i = \{x\in C: x_i=0\}, \forall i\in I_N.
\]

\begin{theorem}\label{The4.1} Assume that the following conditions hold for (\ref{e1}) with $N=3$.
	\begin{itemize}
		\item[(a)] The conditions of Theorem \ref{The2.2} are met so (\ref{e1}) has a unique modified carrying simplex $\Sigma \subset [0, r]$.
		\item[(b)] The three axial fixed points $Q_1, Q_2, Q_3$ are the only fixed points of (\ref{e1}) on $\partial\Sigma = \Sigma\cap \partial C$ and either the inequalities (\ref{e9}) or (\ref{e10}) hold:
		\begin{eqnarray}
			\forall i\in I_3, \forall j\in I_3\setminus\{i, i+1\}, &\;& f_i(Q_{i+1}) < 1 < f_{j}(Q_{i+1}), \label{e9}\\
			\forall i\in I_3, \forall j\in I_3\setminus\{i, i+1\}, & & f_{i+1}(Q_i) < 1 < f_j(Q_i).\label{e10}
		\end{eqnarray}
		\item[(c)] For all $i\in I_3$, $\frac{\partial f_i}{\partial x_i}(Q_i) <0$.
		\item[(d)] There is a unique fixed point $p$ in $\textup{int}\Sigma = \Sigma\setminus \partial\Sigma$ that is hyperbolic with one-dimensional stable manifold $W^s(p)$ in $\dot{C}$ and globally repelling on $\Sigma$.
	\end{itemize}
	Then $\partial\Sigma$ is a globally attracting heteroclinic limit cycle, so $\omega(x^0) = \partial\Sigma$ for all $x^0\in\dot{C}\setminus W^s(p)$.
\end{theorem}

\begin{proof}
	By conditions (a) and (d), we have $\omega(x^0) \subset \partial\Sigma$ for all $x^0\in\dot{C}\setminus W^s(p)$. We need only check that $\partial\Sigma$ is a heteroclinic cycle and prove that $\omega(x^0) = \partial\Sigma$.
	
	Suppose the inequalities (\ref{e9}) hold. From (\ref{e3}), (\ref{e4}) and (\ref{E4}) we know that $DT(Q_1)$ has eigenvalues $f_3(Q_1) <1$, $f_2(Q_1) >1$ and $1+q_1\frac{\partial f_1}{\partial x_1}(Q_1)$. As $\rho(M(Q_1)) = \rho(\tilde{M}(Q_1))=-q_1\frac{\partial f_1}{\partial x_1}(Q_1)$, by (c) and condition (iii) of Theorem \ref{The2.2} we have $0<-q_1\frac{\partial f_1}{\partial x_1}(Q_1)<1$ so $0<1+q_1\frac{\partial f_1}{\partial x_1}(Q_1)<1$. Thus, $Q_1$ is a saddle point with stable manifold $W^s(Q_1)=\pi_2\setminus\pi_1$ and unstable manifold $W^u(Q_1)=\Sigma\cap(\pi_3\setminus\pi_1)$, which connects $Q_1$ to $Q_2$. Similarly, $Q_2$ is a saddle point with stable manifold $W^s(Q_2)=\pi_3\setminus\pi_2$ and unstable manifold $W^u(Q_2)=\Sigma\cap(\pi_1\setminus\pi_2)$, which connects $Q_2$ to $Q_3$; $Q_3$ is a saddle point with stable manifold $W^s(Q_3)=\pi_1\setminus\pi_3$ and unstable manifold $W^u(Q_3)=\Sigma\cap(\pi_2\setminus\pi_3)$, which connects $Q_3$ to $Q_1$. Therefore, $ \partial\Sigma$ is a heteroclinic cycle with the direction $Q_1 \to Q_2\to Q_3\to Q_1$.
	
	Next, we prove that $\omega(x^0) = \partial\Sigma$. Since $\omega(x^0)$ is nonempty, closed, invariant and $\omega(x^0) \subset \partial\Sigma$, there is a $y\in\omega(x^0)\cap\partial\Sigma$. We have either $y\in\{Q_1, Q_2, Q_3\}$ or both $\lim_{n\to\infty}T^n(y)$ and $\lim_{n\to\infty}T^{-n}(y)$ in $\{Q_1, Q_2, Q_3\}$. Thus, $\omega(x^0)$ contains at least one of the $Q_i$. Without loss of generality, we suppose $Q_1\in\omega(x^0)$. That $x^0\not\in W^s(Q_1)$ implies $\omega(x^0)\not= \{Q_1\}$. Since $Q_1$ is a saddle point and $x^0\not\in W^s(Q_1)\cup W^u(Q_1)$, we must have $W^u(Q_3)= W^s(Q_1)\cap \partial\Sigma \subset\omega(x^0)$ and $W^u(Q_1)= W^u(Q_1)\cap \partial\Sigma \subset\omega(x^0)$. So $Q_2, Q_3\in\omega(x^0)$. By the same reasoning as above, we also have $W^u(Q_2)\subset\omega(x^0)$. As
	\[
	\partial\Sigma =\{Q_1, Q_2, Q_3\}\cup W^u(Q_1) \cup W^u(Q_2) \cup W^u(Q_3),
	\]
	we have shown $\partial\Sigma \subset\omega(x^0)$. Hence, $\omega(x^0) = \partial\Sigma$.  
	
	If the inequalities (\ref{e10}) hold, the same reasoning as above is valid with the heteroclinic cycle $\partial\Sigma$ having the direction $Q_1 \to Q_3\to Q_2\to Q_1$ and $\omega(x^0) = \partial\Sigma$.
\end{proof}

\begin{theorem}\label{The4.2} Assume that (a)--(c) in Theorem \ref{The4.1} and the following condition hold for (\ref{e1}) with $N=3$.
	\begin{itemize}
		\item[(d)] There is a unique interior fixed point $p$ that is globally attracting in $\dot{C}$.
	\end{itemize}
	Then $\partial\Sigma$ is a heteroclinic limit cycle globally repelling on $\Sigma$, so $\alpha(x^0) = \Gamma_0$ for all $x^0\in\textup{int}\Sigma\setminus \{p\}$.
\end{theorem}

\begin{proof}
	In the proof of Theorem \ref{The4.1}, replacing $\omega(x^0)$ by $\alpha(x^0)$, the argument is still valid. 
\end{proof}

\section{An example}
In this section, we are going to consider the three-dimensional competitive Ricker model where the map $T$ of (\ref{e1}) is defined by
\begin{equation}\label{e11}
	\forall i\in I_3, T_i(x) = x_i e^{u(1-x_i -\alpha x_{i+1})}, \alpha >1, 0<u<(1+\alpha)^{-1}, x_4=x_1.
\end{equation}
Since $f_i(x) = e^{u(1-x_i -\alpha x_{i+1})}$, we have 
\begin{equation}\label{e12}
	Df(x) = \left(\begin{array}{ccc} \frac{\partial f_1(x)}{\partial x_1} & \frac{\partial f_1(x)}{\partial x_2} & \frac{\partial f_1(x)}{\partial x_3} \\
		\frac{\partial f_2(x)}{\partial x_1} & \frac{\partial f_2(x)}{\partial x_2} & \frac{\partial f_2(x)}{\partial x_3} \\
		\frac{\partial f_3(x)}{\partial x_1} & \frac{\partial f_3(x)}{\partial x_2} & \frac{\partial f_3(x)}{\partial x_3}  \end{array}\right) = \left(\begin{array}{ccc} -uf_1(x) & -u\alpha f_1(x) & 0 \\
		0 & -uf_2(x) & -u\alpha f_2(x) \\
		-u\alpha f_3(x) & 0 & -uf_3(x)  \end{array}\right).
\end{equation}

Note that since 
\[
\frac{\partial f_1(x)}{\partial x_3} = \frac{\partial f_2(x)}{\partial x_1} = \frac{\partial f_3(x)}{\partial x_2} = 0,
\]
we cannot apply Theorem \ref{The1.2} to system (\ref{e1}) with $T$ defined by (\ref{e11}) as its condition (ii) cannot be met. This example shows the situation where Theorem \ref{The2.2} is applicable but Theorem \ref{The1.2} is not. 

If we take $q = (1, 1, 1)$ and any $r\gg q$, the conditions of (i) and (ii) of Theorem \ref{The2.2} are satisfied. For $x\in [0, q]\setminus\{0\}$, it follows from $0<u<(1+\alpha)^{-1}$ and $x_1\leq 1$ that 
\[
f_i(x) +\sum^3_{j=1}x_i\frac{\partial f_i(x)}{\partial x_j} = f_i(x)[1-x_iu(1+\alpha)]>0
\]
for all $i\in I_3$. So (\ref{e7}) holds, which implies $\rho(M(x)) <1$. Hence (iii) of Theorem \ref{The2.2} is met. Therefore, by Theorem \ref{The2.2}, (\ref{e1}) with (\ref{e11}) has a unique modified carrying simplex $\Sigma\subset [0, r]$. 

For $x\in\pi_1$, $f_2(x) =1$ if and only if $x_2+\alpha x_3 =1$ and $f_3(x) =1$ if and only if $x_3=1$. Thus, $Q_2 = (0, 1, 0)$ and $Q_3 = (0, 0, 1)$ are the only two fixed points of $T$ on $\pi_1\setminus\{0\}$. Similarly, $Q_1 = (1, 0, 0)$ and $Q_3$ are the only two fixed points on $\pi_2\setminus\{0\}$, $Q_1$ and $Q_2$ are the only two fixed points on $\pi_3\setminus\{0\}$, so $T$ has only three axial fixed points $Q_1, Q_2, Q_3$ on $\partial\Sigma$. As
\[
f_1(Q_2) = f_2(Q_3) = f_3(Q_1) = e^{u(1-\alpha)}<1,
\]
\[
f_2(Q_1) = f_3(Q_2) = f_1(Q_3) = e^u >1,
\]
the inequalities (\ref{e9}) hold. Thus, the conditions (a)--(c) of Theorem \ref{The4.1} are fulfilled.

The map $T$ has a unique interior fixed point $p = (\frac{1}{1+\alpha}, \frac{1}{1+\alpha}, \frac{1}{1+\alpha})$ with
\[
DT(p) = \left(\begin{array}{ccc} 1-\frac{u}{1+\alpha} & -\frac{u\alpha}{1+\alpha} & 0 \\
	0 & 1-\frac{u}{1+\alpha} & -\frac{u\alpha}{1+\alpha}\\
	-\frac{u\alpha}{1+\alpha} & 0 & 1-\frac{u}{1+\alpha}
\end{array}\right).
\]
The matrix has three eigenvalues:
\begin{equation}\label{e13}
	\lambda_1 = 1-u, \lambda_{2, 3} = 1 + \frac{u(\alpha -2)}{2(1+\alpha)} \pm i \frac{\sqrt{3}u\alpha}{2(1+\alpha)}.
\end{equation}
Under the conditions given in (\ref{e11}), $0<\lambda_1 <1$. To determine whether $|\lambda_{2, 3}|>1$, we define the function
\begin{equation} \label{e14}
	F(\alpha)= \frac{1}{1+\alpha} -\frac{2+\alpha-\alpha^2}{1-\alpha+\alpha^2} = \frac{\alpha^3 + \alpha^2 -4\alpha -1}{(1+\alpha)(1-\alpha+\alpha^2)}, \alpha \geq 0.
\end{equation} 
For system (\ref{e1}) with (\ref{e11}) we have the following result:

\begin{theorem}
	For system (\ref{e1}) with $T$ defined by (\ref{e11}), the following conclusions hold.
	\begin{itemize}
		\item[(a)] There is an $\alpha_0\in (1, 2)$ such that $\alpha_0^3 + \alpha_0^2 -4\alpha_0 -1 = 0$ so $F(\alpha_0)=0$, $F(\alpha) <0$ for $\alpha\in [0, \alpha_0)$ and $F(\alpha) >0$ for $\alpha> \alpha_0$. For $\alpha >\alpha_0$ and $\frac{2+\alpha-\alpha^2}{1-\alpha+\alpha^2} < u < \frac{1}{1+\alpha}$, $|\lambda_{2, 3}|>1$ so $p$ is a repellor on $\Sigma$; for $0 \leq \alpha <\alpha_0$ and $0<u<\frac{1}{1+\alpha}$, $|\lambda_{2, 3}|<1$ so $p$ is asymptotically stable.
		\item[(b)] For $\alpha\geq 2$, $p$ is a global repellor on $\Sigma$.
		\item[(c)] For $\alpha\in [0, \alpha_0)$ with 
		\begin{equation}\label{e18}
			3u(1-\alpha +\alpha^2) < 2+\alpha-\alpha^2,
		\end{equation}	
		$p$ is globally asymptotically stable in $\dot{C}$.
		\item[(d)] For $\alpha \geq 2$ and $0 <u<\frac{1}{1+\alpha}$, $\partial\Sigma$ is a globally attracting heteroclinic limit cycle.
		\item[(e)] For $\alpha\in (1, \alpha_0)$ and $0 < u < \min\{\frac{1}{1+\alpha}, \frac{2+\alpha-\alpha^2}{3(1-\alpha+\alpha^2)}\}$, $\partial\Sigma$ is a heteroclinic limit cycle globally repelling on $\Sigma$.
	\end{itemize}	
\end{theorem}

\begin{proof}
	(a) From (\ref{e13}) we have $|\lambda_{2, 3}| > 1$ if and only if 
	\[
	u > \frac{(2-\alpha)(1+\alpha)}{1-\alpha + \alpha^2} = \frac{2+\alpha-\alpha^2}{1-\alpha+\alpha^2}.
	\]
	This holds obviously if $\alpha \geq 2$ since $0 < u < \frac{1}{1+\alpha}$. For $\alpha\in [0, 2)$, from (\ref{e14}) we see that the denominator of $F(\alpha)$ is positive and the numerator of $F(\alpha)$ has a minimum at $\alpha = \frac{\sqrt{13}-1}{3}$, decreases for $\alpha\in [0, \frac{\sqrt{13}-1}{3})$ and increases for $\alpha > \frac{\sqrt{13}-1}{3}$. As the numerator is negative at $\alpha= 0$ and $\alpha=1$ but positive at $\alpha = 2$, there is an $\alpha_0\in (1, 2)$ such that $\alpha_0^3 + \alpha_0^2 -4\alpha_0 -1 = 0$ so $F(\alpha_0)=0$, $F(\alpha) <0$ for $\alpha\in [0, \alpha_0)$ and $F(\alpha) >0$ for $\alpha> \alpha_0$. Then the conclusion (a) follows.
	
	(d) The conclusion follows from (b) and Theorem \ref{The4.1}.
	
	(e) The conclusion follows from (c) and Theorem \ref{The4.2}.
	
	(b) Let 
	\[
	V(x) = x_1x_2x_3, \quad W(x) = x_1 +x_2 + x_3, \quad x\in C.
	\]
	Then
	\begin{eqnarray}
		V(T(x)) &=& V(x)e^{u(3-(1+\alpha)(x_1+x_2+x_3))} \nonumber \\
		&=& V(x)\exp\{u(1+\alpha)(\frac{3}{1+\alpha}-W(x))\}. \label{e15}
	\end{eqnarray}
	Thus, for all $x\in \dot{C}\setminus\{p\}$, $W(x) > \frac{3}{1+\alpha}$ implies $V(T(x)) < V(x)$, $W(x) < \frac{3}{1+\alpha}$ implies $V(T(x)) > V(x)$, and $W(x) = \frac{3}{1+\alpha}$ implies $V(T(x)) = V(x)$. For all $x\in C$ that is not a fixed point of $T$, 
	\begin{eqnarray*}
		W(T(x)) &=& x_1e^{u(1-x_1-\alpha x_2)} + x_2e^{u(1-x_2-\alpha x_3)} + x_3e^{u(1-x_3-\alpha x_1)}\\
		&> & x_1[1+u(1-x_1-\alpha x_2)] + x_2[1+u(1-x_2-\alpha x_3)] \\
		&& + x_3[1+u(1-x_3-\alpha x_1)] \\
		&=& W(x)(1+u) -u[x_1^2 + x_2^2 + x_3^2 + \alpha x_1x_2 +\alpha x_2x_3 +\alpha x_1x_3] \\
		&=& W(x)(1+u) - uW(x)^2 -u(\alpha -2)(x_1x_2 +x_1x_3 + x_2x_3).
	\end{eqnarray*}	
	As $x_1x_2 +x_1x_3 + x_2x_3 \leq x_1^2 + x_2^2 + x_3^2$, we have
	\begin{eqnarray*}
		x_1x_2 +x_1x_3 + x_2x_3 &\leq& \frac{1}{3}[ 2(x_1x_2 +x_1x_3 + x_2x_3) + x_1^2 + x_2^2 + x_3^2]\\
		&=& \frac{1}{3}W(x)^2.
	\end{eqnarray*}	
	Thus, if $\alpha \geq 2$ then
	\begin{eqnarray}
		W(T(x)) &> & (1+u)W(x) -\frac{u(1+\alpha)}{3} W(x)^2 \nonumber \\
		&=& W(x) + \frac{u(1+\alpha)}{3}W(x)(\frac{3}{1+\alpha} - W(x)).  \label{e16}
	\end{eqnarray}
	Therefore, for all $x\in C$ which is not a fixed point such that $W(x) \leq \frac{3}{1+\alpha}$, we have $W(T(x)) >W(x)$. The conditions of Theorem \ref{The2.2} actually ensure that $T$ from $[0, r]$ to $T([0, r])$ is a homeomorphism (see Remark (2.1) (b) in \cite{Hou3}). So, for each $x\in T([0, r])$, $T^{-1}(x)$ exists in $[0, r]$. Let
	\begin{equation}\label{e17}
		S = \{x\in C:  W(x)\leq \frac{3}{1+\alpha}\}.
	\end{equation}
	By (\ref{e16}) and continuity of $T$ we have $S\subset T(S)$. Since $\overline{\mathcal{B}(0)}= \Sigma\cap\mathcal{B}(0)\subset [0, r]$ is invariant under $T$ so that $T^{-1}(\overline{\mathcal{B}(0)}) = \overline{\mathcal{B}(0)}$, we have
	\[
	T^{-1}(S\cap\overline{\mathcal{B}(0)}) \subset S\cap\overline{\mathcal{B}(0)}.
	\]
	Thus, for each $x\in S\cap\overline{\mathcal{B}(0)}$ which is not a fixed point and for all $n\in\N$, $T^{-n}(x) \in S\cap\overline{\mathcal{B}(0)}$ and $W(T^{-(n+1)}(x))<W(T^{-n}(x))$. Since $S$ contains only two fixed points $0$ and $p$, we must have $\lim_{n\to\infty}W(T^{-n}(x)) =0$ so $\lim_{n\to\infty}T^{-n}(x) =0$. Hence, for all  $x\in S\cap\overline{\mathcal{B}(0)}$, if it is not a fixed point then $x\in \mathcal{B}(0)$. This shows that $(\Sigma\setminus\{p, Q_1, Q_2, Q_3\})\cap S = \emptyset$. Now for each $x\in \textup{int}\Sigma\setminus\{p\}$ and all $n\in\N$, $T^n(x) \in\Sigma\setminus\{p\}$ so $W(T^n(x)) >\frac{3}{1+\alpha}$ and $V(T^{n+1}(x))<V(T^{n}(x))$. Therefore, we must have $\lim_{n\to\infty}V(T^{n}(x)) =0$, so $\omega(x)\subset \partial\Sigma$ and $p$ is a global repellor on $\Sigma$.
	
	(c) Now for $\alpha \in [0, \alpha_0)$ and $x\in C$ not a fixed point, as $0 < u <1$, we have
	\begin{eqnarray*}
		W(T(x)) &=& \frac{x_1}{e^{-u(1-x_1-\alpha x_2)}} + \frac{x_2}{e^{-u(1-x_2-\alpha x_3)}} + \frac{x_3}{e^{-u(1-x_3-\alpha x_1)}}\\
		&< & \frac{x_1}{1-u(1-x_1-\alpha x_2)} + \frac{x_2}{1-u(1-x_2-\alpha x_3)} \\
		&& + \frac{x_3}{1-u(1-x_3-\alpha x_1)} = H(x).
	\end{eqnarray*}
	Define the set $S(\theta) = \{x\in C: W(x)= 3\theta\}$ for each $\theta \in (0, \frac{1}{1+\alpha}]$. It can be shown that $p_0 = (\theta, \theta, \theta)$ is the unique maximum point of $H(x)$ for $x\in S(\theta)$. The proof uses the same technique as that used in \cite{BaHo}. For convenience, we put the proof in the appendix. Clearly, $p_0 = p$ for $\theta = \frac{1}{1+\alpha}$.
	
	Note that $H(p_0)= \frac{3\theta}{c(\theta)}$ as a function of $\theta\in (0, \frac{1}{1+\alpha}]$ is strictly increasing, $H(p_0) < H(p)= \frac{3}{1+\alpha}$ for $\theta < \frac{1}{1+\alpha}$. So far we have proved that $W(T(x)) <H(x) \leq H(p_0)\leq H(p)$ for all $x\in S$ that is not a fixed point for the set $S$ defined by (\ref{e17}). This means that $T(S)\subset S$. 
	
	Now for any $x\in \dot{C}\setminus\{p\}$, we show that $\omega(x) = \{p\}$ so that $p$ is globally attracting in $\dot{C}$. The global asymptotic stability of $p$ follows from this and part (a). If $x\in S$ then $x\in S\setminus(\{p\}\cup \partial C)$, so $T^n(x) \in S\setminus(\{p\}\cup \partial C)$ and $W(T^n(x)) <W(x)\leq \frac{3}{1+\alpha}$. Thus, $V(T^n(x))$ is bounded and increasing for all $n\in\N$. This leads to the existence of a constant $c>0$ such that $\lim_{n\to\infty}V(T^n(x))=c$ and $V(\omega(x)) =c$. As $\omega(x)\subset S\setminus\partial C$, if there is a $p_1\in\omega(x)\setminus\{p\}$ then $p_1$ is not a fixed point, so $W(T(p_1))<\frac{3}{1+\alpha}$ and $V(T^2(p_1)) <V(T(p_1)) \leq V(p_1)=c$. But this contradicts $V(T^2(p_1))=c$ due to $T^2(p_1)\in\omega(x)$ by the invariance of $\omega(x)$. Hence, $\omega(x)\setminus\{p\}= \emptyset$ and $\omega(x)=\{p\}$. If $x\not\in S$ but there is an $m>0$ such that $T^m(x)\in S$, then the above reasoning is valid for $T^m(x)$ so $\omega(x)= \omega(T^m(x))=\{p\}$. If $T^n(x)\not\in S$ for all $n\in\N$, then $V(T^n(x))$ is decreasing so $\lim_{n\to\infty}V(T^n(x)) = d$ for some $d\geq 0$. If $d=0$ then $\omega(x)\subset \partial\Sigma$. By the same reasoning as that used in the proof of Theorem 5 we derive  $\omega(x)= \partial\Sigma$ based on the fact that $\partial\Sigma$ is a heteroclinic cycle. In particular, $Q_1 = (1, 0, 0) \in \omega(x)$. Note that $W(Q_1) = 1 < \frac{3}{1+\alpha}$. This contradicts the fact that $W(y) \geq \frac{3}{1+\alpha}$ for all $y\in\omega(x)$ due to $W(T^n(x)) >\frac{3}{1+\alpha}$. This contradiction shows that $d>0$. Then, by the same argument as above, we must have $\omega(x) = \{p\}$. The proof of the theorem is complete.
\end{proof}

\section{Conclusion}
We have reviewed the modified carrying simplex theory as a recent development in this area. By comparing the concept of modified carry simplex and the sufficient conditions for system (\ref{e1}) with (\ref{e2}) to have a modified carrying simplex with those of carrying simplex, we note that the modified carrying simplex keeps most of the features for carrying simplex but requires much weaker conditions on the map $T$: it requires $\frac{\partial f_i}{\partial x_j} \leq 0$ rather than $\frac{\partial f_i}{\partial x_j} < 0$ for all $i, j\in I_N$ and all $x\in [0, r]$. This means that, like the theorems for carrying simplex, the theorems for modified carrying simplex can be applied to a broader class of systems in various applications.

To demonstrate this point, we have found a criterion for a three dimensional system to have a globally attracting (repelling) heteroclinic limit cycle when it admits a modified carrying simplex. Then we have shown by a Ricker model as a concrete example that the theorem for modified carrying simplex is applicable but the theorem for carrying simplex is not.

We expect that various applications of the theorems for modified carrying simplex will appear in future.

\section*{Appendix: Proof that $H(x) < H(p_0)$ for $x\in S(\theta)\setminus\{p_0\}$}

For this purpose, we first show that 
\[
H(x)< H(p_0) = \frac{3\theta}{1-u+u(1+\alpha)\theta}=\frac{3\theta}{c(\theta)}, \quad x\in \partial S(\theta) = S(\theta)\cap \partial C. 
\]
As $H(x)$ has the rotational symmetry about the components of $x$, we need only show the inequality for $x\in S(\theta)$ with $x_3 = 0$ and $x_1 + x_2 = 3\theta$. Clearly, for $x =(3\theta, 0, 0)$ we have $H(x) = \frac{3\theta}{1+u(3\theta -1)} < H(p_0)$ due to $ 3 > 1 + \alpha$. If $x_3 = 0$ but $x_2\not= 0$, then  	
\begin{eqnarray*}
	H(x)-H(p_0) &=& \frac{x_1}{1-u(1-x_1-\alpha x_2)} + \frac{x_2}{1-u(1-x_2)} -\frac{3\theta}{c(\theta)} \\
	&=& H(p_0)\left\{\frac{c(\theta)x_1}{(1-u)(x_1+x_2)+3u\theta(x_1+\alpha x_2)}\right. \\
	&& \left.  + \frac{c(\theta)x_2}{(1-u)(x_1+x_2)+ 3u\theta x_2)} -1 \right\} \\
	&=& H(p_0)\left\{\frac{c(\theta)(x_1/x_2)}{(1+u(3\theta-1))(x_1/x_2)+ (1+u(3\theta\alpha -1))}\right. \\
	&& \left.  + \frac{c(\theta)}{(1-u)(x_1/x_2)+ (1+u(3\theta-1))} -1 \right\} \\
	&=& H(p_0)\left\{\frac{c(\theta)X}{AX+B} + \frac{c(\theta)}{CX+A} -1\right\}\\
	&=& \frac{H(p_0)R(X)}{(AX+B)(CX+A)},
\end{eqnarray*}
where $A = 1-u+3u\theta >c(\theta)>0$, $B= 1-u+3u\theta\alpha >0$, $C=1-u>0$, $X = x_1/x_2 \geq 0$ and 
\[
R(X) = C(c(\theta)-A)X^2-(A^2-2c(\theta)A+BC)X+B(c(\theta)-A).
\]
As $C(c(\theta)-A) <0$ and $B(c(\theta)-A)<0$, we have $H(x) - H(p_0) <0$ for all $x\in \partial S(\theta)$ if either $A^2-2c(\theta)A+BC \geq 0$ or
\begin{eqnarray*}
	[A^2-2c(\theta)A+BC]^2 &-& 4BC(c(\theta)-A)^2 \\
	&=& (A^2-BC)[A^2-BC-4c(\theta)(A-c(\theta))] <0.
\end{eqnarray*}
The former is reduced to $3u\theta(1-2\alpha)\geq (1-u)(2-\alpha)>0$, which is impossible for $\alpha\geq 0.5$ or small $u\theta$. Since
\begin{eqnarray*}
	A^2-BC &=& (1-u+3u\theta)^2-(1-u+3u\theta\alpha)(1-u) \\
	&=& 3u(1-u)\theta(2-\alpha) + 9\theta^2u^2 >0,
\end{eqnarray*}
the latter becomes $A^2-BC-4c(\theta)(A-c(\theta))<0$, which is reduced to
\begin{equation}\label{e19}
	u\theta(2\alpha-1)^2 < (2-\alpha)(1-u).	
\end{equation}
When $\theta$ is replaced by $\frac{1}{1+\alpha}$, (\ref{e19}) is equivalent to (\ref{e18}). As $0 <\theta \leq \frac{1}{1+\alpha}$, (\ref{e19}) holds under the condition (\ref{e18}). Therefore, $H(x)<H(p_0)$ for all $x\in\partial S(\theta)$.

Next, we show that $H(x) < H(p_0)$ for all $x\in S(\theta)\setminus\{p_0\}$. Clearly, $H(x)-H(p_0) =0$ for $x=p_0$. For any $x_0\in \textup{int} S(\theta)\setminus\{p_0\}$, $x_0$ and $p_0$ determine a unique line segment $\overline{yz}$ with $y, z\in\partial S(\theta)$ such that $x_0, p_0 \in \overline{yz}$. The points on $\overline{yz}$ can be written as $x(s) = y + s(z-y)$ for $s\in [0, 1]$ such that $x(0) = y, x(1) = z$, $x(s_1) = x_0$ and $x(s_2) = p_0$ for some distinct $s_1, s_2 \in (0, 1)$. Let $h(s) = H(x(s)) -H(p_0)$. From the definition of $H(x)$ we can write it as $H(x) = \frac{P(x)}{Q(x)}$, where
\begin{eqnarray*}
	P(x) &=& x_1[1-u(1-x_2-\alpha x_3)][1-u(1-x_3-\alpha x_1)] \\
	&& +x_2[1-u(1-x_1-\alpha x_2)][1-u(1-x_3-\alpha x_1)] \\
	&& + x_3[1-u(1-x_2-\alpha x_3)][1-u(1-x_1-\alpha x_2)], \\
	Q(x) &=& [1-u(1-x_2-\alpha x_3)][1-u(1-x_2-\alpha x_3)][1-u(1-x_3-\alpha x_1)].
\end{eqnarray*}
Then $h(s) = \frac{P(x(s)) -H(p_0) Q(x(s))}{Q(x(s))}= \frac{h_1(s)}{Q(x(s))}$ with $Q(x) >0$. As $h(0) = H(y)-H(p_0) <0$, $h(1) = H(z)-H(p_0) <0$ and $h(s_2) = H(p_0)-H(p_0) =0$, the polynomial $h_1(s)$ of degree three satisfies $h_1(0)<0$, $h_1(1) <0$ and $h_1(s_2) = 0$. Also, from
\[
\frac{d P(x(s))}{ds}|_{s=s_2} = [c(\theta)^2+2u\theta c(\theta)(1+\alpha)] (W(z)-W(y)) =0
\]
and
\[
\frac{d Q(x(s))}{ds}|_{s=s_2} = uc(\theta)^2(1+\alpha) (W(z)-W(y)) =0
\]
we obtain $h_1'(s_2) = 0$. Thus, $h_1(s) = (s-s_2)^2(Ds+E)$. As $h_1(0)= (s_2)^2E<0$ and $h_1(1) = (1-s_2)^2(D+E)<0$, these imply that $Ds+E<0$ for $s\in [0, 1]$, so $h_1(s_1)<0$ and $h(s_1) = H(x_0)-H(p_0)<0$. Hence, $H(x) < H(p_0)$ for all $x\in S(\theta)\setminus\{p_0\}$.

\section*{Acknowledgement} The author thanks the referee for providing positive comments and suggestions adopted in this version of the paper.
%
%


\begin{thebibliography}{99}
	%
	\bibitem{Bai1} Baigent, S.: \emph{Convexity of the carrying simplex for discrete-time planar competitive Kolmogorov systems}, J. Differ. Equ. Appl., \textbf{22 (5)} (2016), 609--620.
	
	\bibitem{Bai2} Baigent, S.: \emph{Convex geometry of the carrying simple for the May-Leonard map}, Discret. Contin. Dyn. Syst. B, \textbf{24 (4)} (2019), 1697--1723.
	
	\bibitem{BaHo} Baigent, S. and Hou, Z.: \emph{Global stability of discrete-time competitive population models}, J. Difference Equ. Appl., \textbf{23} (2017), 1378--1396.
	
	\bibitem{BaElLu} Cabral Balreira, E., Elaydi, S. and Lu\'{i}s, R.: \emph{Global stability of higher dimensional monotone maps}, J. Difference Equ. Appl., \textbf{23} (2017), 2037--2071.
	
	\bibitem{DiWaYa} Diekmann, O., Wang, Y. and Yan, P.: \emph{Carrying simplices in discrete competitive systems and age-structured semelparous populations}, Discrete Contin. Dyn. Syst., \textbf{20} (2008), 37--52.
	
	\bibitem{GyJiNiYa2} Gyllenberg, M., Jiang, J., Niu, L. and Yan, P.: \emph{On the classification of generalised competitive Atkinson-Allen models via the dynamics on the boundary of the carrying simplex}, Discret. Contin. Dyn. Syst., \textbf{38} (2018), 615--650.
	
	\bibitem{GyJiNi3} Gyllenberg, M., Jiang, J. and Niu, L.: \emph{A note on global stability of three-dimensional Ricker models}, J. Difference Equ. Appl., \textbf{25} (2019), 142--150.
	
	\bibitem{GyJiNiYa1} Gyllenberg, M., Jiang, J., Niu, L. and Yan, P.:
	\emph{On the dynamics of multi-species Ricker models admitting a carrying simplex}, J. Difference Equ. Appl., \textbf{25} (2019), 1489--1530.
	
	\bibitem{Her} Ruiz-Herrera, A.: \emph{Exclusion and dominance in discrete population models via the carrying simplex}, J. Differ. Equ. Appl. \textbf{19} (1) (2013), 96--113.
	
	\bibitem{Hir1} Hirsch, M. W.: \emph{Systems of differential equations that are competitive or cooperative. III: Competing species}, Nonlinearity \textbf{1} (1988), 51-71.
	
	\bibitem{Hir2} Hirsch, M. W.: \emph{On existence and uniqueness of the carrying simplex for competitive dynamical systems}, J. Biol. Dyn. \textbf{2 (2)} (2008), 169--179.
	
	\bibitem{Hou} Hou, Z.: \emph{Geometric method for global stability of discrete population models}, Discret. Contin. Dyn. Syst. B, \textbf{25 (9)} (2020) 3305--3334.
	
	\bibitem{Hou1} Hou, Z.: \emph{On existence and uniqueness of a carrying simplex in Kolmogorov differential systems}, Nonlinearity \textbf{33} (2020) 7067--7087.
	(https://doi.org/10.1088/1361-6544/abb03c)
	
	\bibitem{Hou2} Hou, Z.: \emph{Corrigendum: On existence and uniqueness of a carrying simplex in Kolmogorov differential systems (2020 Nonlinearity 33 7067)}, Nonlinearity \textbf{34} (2021) C5--C6.
	(https://doi.org/10.1088/1361-6544/ac2a4f)
	
	\bibitem{Hou3} Hou, Z.: \emph{On existence and uniqueness of a modified carrying simplex for discrete Kolmogorov systems}, J. Difference Equ. Appl., \textbf{27 (2)} (2021) 284-315.
	(https://doi.org/10.1080/10236198.2021.1894141)
	
	\bibitem{JiMiWa} Jiang, J., Mierczy\'{n}ski, J. and Wang, Y.: \emph{Smoothness of the carrying simplex for discrete-time competitive dynamical systems: A characterization of neat embedding}, J. Differ. Equ. \textbf{246} (4) (2009), 1623--1672.
	
	\bibitem{JiNi1} Jiang, J. and Niu, L.: \emph{On the equivalent classification of three-dimensional competitive Atkinson-Allen models relative to the boundary fixed points}, Discret. Contin. Dyn. Syst. \textbf{36} (1) (2016), 217--244.
	
	\bibitem{JiNi2} Jiang, J. and Niu, L.: \emph{On the equivalent classification of three-dimensional competitive Leslie-Gower models via the boundary dynamics on the carrying simplex}, J. Math. Biol. \textbf{74} (2017), 1223--1261.
	
	\bibitem{JiNiWa} Jiang, J., Niu, L. and Wang, Y.: \emph{On heteroclinic cycles of competitive maps via carrying simplices}, J. Math. Biol. \textbf{72} (2016), 939--972.
	
	\bibitem{WaJi1} Wang, Y. and Jiang, J.: \emph{The general properties of discrete-time competitive dynamical systems}, J. Differ. Equ. \textbf{176} (2001), 470--493.
	
	\bibitem{WaJi2} Wang, Y. and Jiang, J.: \emph{Uniqueness and attractivity of the carrying simplex for discrete-time competitive dynamical systems}, J. Differ. Equ. \textbf{186} (2002), 611--632.
	
	
\end{thebibliography}
\end{document}